\documentclass[]{article}
\usepackage{amsfonts}
\usepackage{dsfont}
\usepackage{amsmath}
\usepackage{verbatim}
\usepackage{amsthm}
\usepackage{enumerate}
\usepackage{bm}
\usepackage{tikz}

\theoremstyle{plain}
\newtheorem{definition}{Definition}
\newtheorem{theorem}{Theorem}
\newtheorem{corollary}{Corollary}
\newtheorem{lemma}{Lemma}

\theoremstyle{definition}
\newtheorem{remark}{Remark}
\newtheorem{example}{Example}

\title{Limits of random trees II\thanks{MSC2010 Subject Classification: 05C80}}
\author{Attila De\'ak\thanks{MTA-ELTE "Numerical Analysis and Large Networks" Research Group}\\
\texttt{deak51@cs.elte.hu}}
\begin{document}
\maketitle

\begin{abstract}
Local convergence of bounded degree graphs was introduced by Benjamini and Schramm.
This result was extended further by Lyons to bounded average degree graphs.
In this paper we study the convergence of random tree sequences with given degree distributions.
Denote by ${\cal D}_n$ the set of possible degree sequences of a labeled tree on $n$ nodes.
Let ${\bm D}_n$ be a random variable on ${\cal D}_n$ and ${\bm T}({\bm D}_n)$ be a uniform random labeled tree with degree sequence ${\bm D}_n$.
We show that the sequence ${\bm T}({\bm D}_n)$ converges in probability if and only if ${\bm D}_n\to {\bm D}=({\bm D}(i))_{i=1}^\infty$, where ${\bm D}(i)\sim {\bm D}(j)$, $\mathds{E}({\bm D}(1))=2$ and ${\bm D}(1)$ is a random variable on $\mathds{N}^+$.\\

\end{abstract}

\section{Introduction}
\label{sec.int}
In recent years the study of the structure and behavior of real world networks has received wide attention.
The degree sequence of these networks appear to have special properties (like power law degree distribution).
Classical random graph models (like the Erd\H{o}s-R\'enyi model) have very different degree sequence.
An obvious solution is to study a random graph with given degree sequence.
More generally generate a random graph with a degree sequence from a family of degree sequences.
In \cite{Chatterjee-Diaconis-Sly} Chatterjee, Diaconis and Sly studied random dense graphs (graphs whose number of edges is comparable to the square of the number of vertices) with a given degree sequence.

It is not always easy to generate a truly random graph with a given degree sequence.
There is a fairly large literature on the configuration model (for the exact definition of the model see \cite{Bollobas}), where for a given degree sequence for each node $i$ we consider $d_i$ stubs and take a random pairing of the stubs and connect the corresponding nodes with an edge.
This model creates the required degree distribution, but gives a graph with possible loops and parallel edges.

A notion of convergence for (dense) graph sequences was developed by Borgs, Chayes, Lov\'asz, S\'os and Vesztergombi in \cite{Borgs-Chayes-Lovasz-Sos-Vesztergombi}.
The limit objects were described by Lov\'asz and Szegedy in \cite{Lovasz-Szegedy}.
Using this limit theory, the authors in \cite{Chatterjee-Diaconis-Sly} described the structure of random (dense) graphs from the configuration model.
They defined the convergence of degree sequences and
for convergent degree sequences they gave a sufficient condition on the degree sequence, which implies the convergence of the random graph sequence (taken from the configuration model).

What can we say if the graphs we want to study are sparse (the number of edges is comparable to the number of vertices) and not dense?
Is there a similar characterization for sparse graphs with given degree sequence?
We establish a characterization for random trees with given (possibly random) degree sequence.
There are various limit theories and convergence notions for trees introduced
by Aldous \cite{Aldous} and by Elek and Tardos \cite{Elek-Tardos}.
We use the notion of convergence introduced for bounded degree graphs (that is the degree of each vertex is bounded above by some uniform constant $d$) first introduced by Benjamini and Schramm \cite{Benjamini-Schramm}.
This notion was extended by Lyons \cite{Lyons} to bounded average degree graphs.

In \cite{Deak} the author described the behavior of a random tree sequence with a given degree distribution.
In this paper we extend this result and prove a similar characterization as in \cite{Chatterjee-Diaconis-Sly} for random trees with given degree sequence.
We define the convergence of degree sequences and give a necessary and sufficient condition on the degree sequence, which implies the convergence of the tree sequence ${\bm T}({\bm D}_n)$ in the sense of Lyons \cite{Lyons}.
In the case of convergence we describe the limit object.

This paper is organized as follows:
In Section \ref{sec.def.not}, we give the basic definitions and notations.
In Section \ref{sec.lim.deg.seq}, we describe the basic properties and the limit of a sequence of random degree sequences.
At the end of the section we state our main theorem.
In Section \ref{sec.lab.subg.dens}, we deal with labeled homomorphisms and in Section \ref{sec.limit}, we describe the limit object.


\section{Basic definitions and notations}
\label{sec.def.not}

\subsection{Random weak limit of graph sequences}

Let $G=G(V,E)$ be a finite simple graph on $n$ nodes.
For $S\subseteq V(G)$ denote by $G[S]$ the subgraph of $G$ spanned by the vertices $v\in S$.
For a finite simple graph $G$ on $n$ nodes , let $B_G(v,R)$
be the rooted $R$-ball around the node $v$, also called as the $R$-neighborhood of $v$, that is the subgraph induced by the nodes at distance at most $R$ from $v$:
$$
B_G(v,R)=G[\{u\in V(G): dist_G(u,v)\leq R\}].
$$
Two rooted graphs $G_1, G_2$ are rooted isomorphic if there is an isomorphism between them which maps the root of $G_1$ to the root of $G_2$.
Given a positive integer $R$, a finite rooted graph $F$ and a probability distribution $\rho$ on rooted graphs, let $p(R,F,\rho)$ denote the probability that the graph $F$ is rooted isomorphic to the $R$-ball around the root of a rooted graph chosen with distribution $\rho$.
It is clear that $p(R,F,\rho)$ depends only on the component of the root of the graph chosen from $\rho$.
So we will assume that $\rho$ is concentrated on connected graphs.
For a finite graph $G$, let $U(G)$ denote the distribution graphs.ed graphs obtained by choosing a uniform random vertex of $G$ as root of $G$.
It is easy to see, that for any finite graph $G$ we have
$$
p(R,F,U(G))={|\{v\in V(G):B_G(v,R)\textrm{ is rooted isomorhpic to } F\}|\over |V(G)|}.
$$

\begin{definition}
Let ($G_n$) be a sequence of finite graphs on $n$ nodes, $\rho$ a probability distribution on infinite rooted graphs. We say that the random weak limit of $G_n$ is $\rho$, if for any positive integer $R$ and finite rooted graph $F$, we have
\begin{equation}
\label{eqn.conv}
\lim_{n\rightarrow \infty}p(R,F,U(G_n))=p(R,F,\rho).
\end{equation}
\end{definition}
If $G_n$ is a sequence of random finite graphs, then $p(R,F,U(G_n))$ is a random variable, so by convergence we mean convergence in probability.
\begin{definition}
Let $(G_n)$ be a sequence of random finite graphs on $n$ nodes, $\rho$ a probability distribution on infinite rooted graphs.
We say that the random weak limit of $G_n$ is $\rho$, if $\forall\epsilon > 0,R\in \mathds{N}^+$ and finite rooted graph $F$, we have
\begin{equation}
\label{eqn.rnd.conv}
\lim_{n\rightarrow \infty}\mathds{P}(|p(R,F,U(G_n))-p(R,F,\rho)|>\epsilon)=0.
\end{equation}
\end{definition}

The formal meaning of this formula is, that the statistics $p(R,F,U(G_n))$ as random variables are concentrated.

\subsection{Other notations}

We will denote random variables with bold characters.
For a probability space $(\Omega,{\cal B},\mu)$ and $A\in {\cal B}$ denote by ${\bm I}(A)$ the indicator variable of the event $A$.
Denote by ${\cal D}_n$ the set of possible degree sequences of a labeled tree on $n$ nodes.
Throughout the paper we consider labeled trees on $n$ nodes unless stated otherwise.
Let ${\bm D}_n$ be a random variable on ${\cal D}_n$.
Denote by ${\bm T}({\bm D}_n)$ the uniform random tree on $n$ labeled nodes, with degree sequence ${\bm D}_n$.
Denote the degree sequence of a tree $T$ by $D_T=(D_T(i))_{i=1}^n$.
For a given degree sequence $D=(D(i))_{i=1}^n$ there are
$$n-2 \choose D(1)-1,D(2)-1,\cdots ,D(n)-1$$
labeled trees with degree sequence $D$.
It follows that for an arbitrary tree $T$
\begin{equation*}
\mathds{P}({\bm T}({\bm D}_n)=T)=\frac{\mathds{P}({\bm D}_n=D_T)}{
\displaystyle{{n-2\choose D_T(1)-1,D_T(2)-1,\cdots ,D_T(n)-1}}}.
\end{equation*}
If it does not cause any confusion, we will use $D,{\bm T}_n$ instead of $D_T,{\bm T}({\bm D}_n)$ respectively.

A finite rooted graph $G$ with root $v$ is said to be $l$-deep if the largest distance from the root is $l$, that is
$$l=\max_{u\in V(G)}dist(v,u).$$

Denote by $U^l$ the set of equivalence classes of finite unlabeled $l$-deep rooted graphs with respect to root-preserving isomorphisms.
Let $T_x^l$ be an $l$-deep rooted tree on $k$ nodes with root $x$.
Denote the vertices at distance $i$ from the root by $T_i$, and let $t_i=|T_i|$ ($t_0$ is $1$, $t_1$ is the degree of the root).
For every finite graph $G$, $p(R,F,U(G))$ induces a probability measure on $U^R$ which we call the $R$-neighborhood statistics of $G$.
If $G$ is a tree then $p(R,F,U(G))$ is concentrated on rooted trees.

Let $\cal T$ be the set of all countable, connected infinite rooted trees.
For an infinite rooted tree $T\in {\cal T}$ denote by $T(R)$ the $R$-neighborhood of the root of $T$.
For an $R$-deep rooted tree $F$ define the set
$$
{\cal T}(F)=\{T\in {\cal T}: T(R) \textrm{ is rooted isomorphic to } F\}.
$$

Let ${\cal F}$ be the sigma-algebra generated by the sets $({\cal T}(F))_F$, where $F$ is an arbitrary finite rooted tree.
$({\cal T},{\cal F})$ is a probability field.
We call a probability measure $\mu$ on ${\cal T}$ an infinite rooted random tree.

Every infinite random tree $\mu$ has the property that for any $F\in U^R$:

\begin{equation}
\label{eqn.neigh.cons}
p(R,F,\mu)=\sum_{H\in U^{R+1},\ H(R)\cong F}p(R+1,H,\mu).
\end{equation}

Actually every distribution on rooted infinite graphs has the above property.
Note that if we want to prove the convergence of a random tree sequence to a certain limit distribution $\rho$, then we need to have the convergence of the neighborhood densities and also (\ref{eqn.neigh.cons}), the consistency of these densities, which ensures that $\rho$ will be concentrated on infinite rooted trees.
These together will imply (\ref{eqn.rnd.conv}).

\section{Limits of degree sequences}
\label{sec.lim.deg.seq}

Consider a random degree sequence ${\bm D}_n=({\bm D}_{n}(i))_{i=1}^n$ and construct a labeled tree ${\bm T}({\bm D}_n)$ with uniform distribution given the degree sequence.
We want to describe the limit of ${\bm T}({\bm D}_n)$ as $n\to \infty$.
We give a characterization of the degree sequences for which ${\bm T}({\bm D}_n)$ has an infinite random tree as a limit.
To describe the model and the limit, we need to define and understand the limit of a random degree sequence ${\bm D}_n$.
Here we only deal with degree sequences of trees.
We further assume that ${\bm D}_n$ is an exchangeable sequence, that is for any $\sigma\in S_n$ we have
$$
({\bm D}_n(i))_{i=1}^n\sim ({\bm D}_n(\sigma(i)))_{i=1}^n.
$$
Exchangeability is a way to eliminate exceptional vertices and allows us to use the limit theory of exchangeable sequences.
For more on exchangeable random variables we refer to \cite{Aldous.ex}.
\begin{definition}
We say that an exchangeable sequence ${\bm D}_n$ is convergent and ${\bm D}_n\to \bm{D}$, where ${\bm D}$ is a random infinite sequence, if for every $k\in \mathds{N}$ we have
\begin{equation*}
({\bm D}_{n}(i))_{i=1}^k\stackrel{\mathds{P}}{\rightarrow} ({\bm D}(i))_{i=1}^k.
\end{equation*}
\end{definition}
It is easy to see, that if ${\bm D}_n$ is exchangeable and ${\bm D}_n\to {\bm D}$ then ${\bm D}$ is also an exchangeable sequence.
The following theorem of Hewitt and Savage( see \cite{Hewitt-Savage}), which is a generalization of de Finetti's theorem, describes the limits of exchangeable sequences.

\begin{theorem}
\label{thm.deFinetti}
Let $\bm X$ be a random infinite exchangeable array.
Then $\bm X$ is a mixture of infinite dimensional iid distributions
\begin{equation*}
{\bm X} =\int_{IID}\lambda dp(\lambda),
\end{equation*}
\noindent
where $p$ is a distribution on infinite dimensional IID distributions $\lambda$.
\end{theorem}
As a result we have that the limit of an exchangeable degree sequence is an infinite exchangeable sequence and so a mixture of IID distributions.
Note that if ${\bm D}_n$ is not exchangeable then we can take a random permutation $\sigma\in S_n$ and define the exchangeable degree sequence $\tilde{{\bm D}}_n(i)={\bm D}_n(\sigma(i))$.

\begin{lemma}
\label{lem.exchange.degen}
Let $\bm X$ be an infinite exchangeable random sequence. Further assume that we have
$$
\mathds{P}({\bm X}(1)=i,{\bm X}(2)=i)=\mathds{P}({\bm X}(1)=i)\mathds{P}({\bm X}(2)=i).
$$
Then $\bm X$ is an infinite IID distribution ($p$ is concentrated on one distribution).
\end{lemma}
\noindent
{\bf Proof:}
From Jensen's inequality we have that
\begin{equation}
\label{eqn.jensen}
\int_{IID}\lambda(i)^2dp(\lambda)\geq \left(\int_{IID}\lambda(i)dp(\lambda)\right)^2.
\end{equation}
Also from Theorem \ref{thm.deFinetti} we have
\begin{multline*}
\label{eqn.exchange.degen}
\int_{IID}\lambda(i)^2dp(\lambda)=\mathds{P}({\bm X}(1)=i,{\bm X}(2)=i)=\\
=\mathds{P}({\bm X}(1)=i)\mathds{P}({\bm X}(2)=i) =\left(\int_{IID}\lambda(i)dp(\lambda)\right)^2
\end{multline*}
It follows that in (\ref{eqn.jensen}) equality holds which means that $p$ is a degenerate distribution and so proves our lemma.

\qed

We will see that if ${\bm T}({\bm D}_n)$ is convergent, then ${\bm D}_n$ satisfies the assumptions in Lemma \ref{lem.exchange.degen}.
So for a convergent random tree sequence ${\bm T}({\cal D}_n)$ the limit of the degree sequence ${\bm D}_n$ needs to be an infinite IID distribution.

\begin{example}
\label{example.star}
Let ${\bm X}$ be a uniform random element of $[n]$.
Consider the degree sequence
\begin{equation*}
{\bm D}(i)=\left \{
\begin{array}{ll}
n-1, & \textrm{if } i={\bm X}\\
1, & \textrm{otherwise}
\end{array}
\right.
\end{equation*}
Let ${\bm S}{\bm t}_n={\bm T}({\bm D}_n)$ be the star-graph on $n$ nodes.
The limit degree sequence is just the constant $1$ vector $\mathds{1}=(1,1,\cdots )$.
Obviously in the limit the expected degree of a node is $1$.
It is not hard to see, that if $F$ is not a single edge, then $u(R,F,{\bm S}{\bm t}_n)=0$ for every $n>|V(F)|$.
Thus there is no limit distribution $\rho$ on infinite graphs such that $\mathds{P}(|p(R,F,U({\bm S}{\bm t}_n)-p(R,F,\rho)|>\epsilon)\to 0$ for every $F$.
\end{example}
Example \ref{example.star} shows that if only the average degree is bounded, too many unbounded degree vertices destroy convergence.
As the average degree of a tree on $n$ nodes is $2{n-1\over n}$, one would expect that in the limit distribution the expected degree of a node is $2$, that is $\mathds{E}({\bm D}(i))=2$ for every $i$.

It turns out that it is enough to have that the degree sequence converges and $\mathds{E}({\bm D}(i))=2$ holds $\forall i$.
Now we are ready to state our main theorem which describes the degree sequence of convergent random tree sequences.

\begin{theorem}
\label{thm.main}
Let ${\bm D}_n$ be a sequence of random degree sequences (${\bm D}_n\in {\cal D}_n$).
The random tree sequence ${\bm T}({\bm D}_n)$ is convergent and converges to an infinite random tree if and only if ${\bm D}_n\to {\bm D}$, where ${\bm D}=({\bm D}_0,{\bm D}_0,\cdots)$ is an infinite IID sequence and $\mathds{E}({\bm D}_0)=2$.\\

\end{theorem}

\section{Labeled subgraph densities}
\label{sec.lab.subg.dens}

To prove convergence we need to understand the neighborhood statistics of the random tree ${\bm T}({\bm D}_n)$.
First we will count subgraph densities and then relate them to neighborhood statistics.
For fixed unlabeled graphs $F$ and $G$ denote by
$$
inj(F,G)={|\{\phi: \phi \textrm{ is an injective homomorphism from } F \textrm{ to } G\}|\over |V(G)|}
$$
the normalized number of copies of $F$ in $G$.
We call $F$ the test graph.
We call $inj(F,G)$ the injective density of $F$ in $G$.
For bounded degree graphs the convergence of injective densities for every $F$ is equivalent to the convergence of neighborhood densities for every $H$ rooted finite graph.
For bounded average degree graphs subgraph statistics may be unbounded.
For the random star tree ${\bm S}{\bm t}_n$ we have
$$
inj(
\begin{tikzpicture}[every node/.style={circle, fill=black, inner sep=0mm,
minimum size=1mm}]
\node (A) at (-0.2,0) {};
\node (B) at (0.2,0) {};
\node (C) at (0,0.2) {};
\draw (A) -- (C) -- (B);
\end{tikzpicture},
{\bm S}{\bm t}_n)={(n-1)(n-2)\over n}.
$$
To avoid unbounded subgraph statistics we add a further structure to the test graph $F$.
We call a pair $(F,r)$ a numbered graph, where $r=(r_i)_{i=1}^{V(F)}$ and $r_i\in \mathds{N}$.
We call $r_i$ the remainder degree of the node $i\in V(F)$.
Let $(F,r)$ be a numbered graph and $\phi$ be a homomorphism from $F$ to a graph $G$.
We say that $\phi$ is a labeled homomorphism if $\phi$ is a homomorphism and
$$D_G(\phi(v))=D_F(v)+r_v,\ \forall v\in V(F).$$
Let
$$
inj_{lab}((F,r),G)={|\{\phi: \phi \textrm{ is an injective labeled homomorphism from } F \textrm{ to } G\}|\over |V(G)|}
$$
be the normalized number of numbered copies of $F$ in $G$.
First we want to derive properties of degree sequences ${\rm D}_n$ for which $inj_{lab}((F,r),{\bm T}({\bm D}_n))$ is convergent for every finite graph $F$ and remainder degrees $r$.
Then in Section \ref{sec.limit} we will turn to the convergence of neighborhood statistics.

\begin{remark}
The convergence of $inj_{lab}(.,G_n)$ for every $(F,r)$ does not imply the random weak convergence of $G_n$ in general.
$inj_{lab}((F,R),{\bm S}{\bm t}_n)$ is convergent for every $(F,r)$, but as we saw earlier ${\bm S}{\bm t}_n$ is not a convergent tree sequence.
\end{remark}

\begin{remark}
\label{rem.lab.hom.bounded}
Let $(F,r)$ be an arbitrary numbered graph on $k$ nodes.
One can easily see that $inj_{lab}((F,r),G)$ is uniformly bounded for every $G$.
\end{remark}
\medskip
\noindent
{\bf Proof.:}
To see this we will bound the number of ways we can construct an injective labeled homomorphism $\psi$ from $(F,r)$ to $G$.
Let $R=\max\{r_i\}$.
If we define $\psi(1)=v\in V(G)$, then $D_G(\psi(1))=r_1$.
There are at most $D_G(\psi(1))^{D_F(1)}=r_1^{D_F(1)}\leq R^k$ possibilities for $\psi(u)$'s $(u\in N_F(1))$, where $N_F(1)$ is the set of neighbors of $1$ in $F$.
Following this idea we get that for every $v$ there are at most $(R^k)^k$ possible ways to extend $\psi$, given $\psi(1)=v$.
Hence there are at most $nR^{k^2}$ injective labeled homomorphisms from $F$ to $G$ and the remark follows.
\qed

\medskip
\noindent
For an arbitrary numbered tree $(T,r)$, and $\phi:V(T)\mapsto [n]$ let
$$
{\bm I}_n((T,r),\phi)={\bm I}(\{\phi \textrm{ is an injective labeled homomorphism of } T \textrm{ to } {\bm T}_n\})
$$
\begin{equation}
\label{eqn.sum.of.indicator}
{\bm X}_n^{(T,r)}=\sum_{\phi: V(T)\mapsto [n]}{\bm I}_n((T,r),\phi)=n\cdot inj_{lab}((T,r),{\bm T}({\bm D}_n)).
\end{equation}
We define ${\bm I}_n((F,r),\phi),{\bm X}^{(F,r)}_n$ similarly for a numbered forest $(F,r)$.
If it does not cause any confusion, we will omit $r$ from the formulas above and use ${\bm I}_n(T,\phi)$,  $X^T_n,$ ${\bm I}_n(F,\phi)$ and ${\bm X}^{F}_n$ instead to simplify notation.
For random graph sequences ${\bm G}_n$ by the convergence of $inj_{lab}(.,{\bm G}_n)$ we mean convergence in probability.

Let ${\bm D}_n$ be a random degree sequence and ${\bm T}_n={\bm T}({\bm D}_n)$ be the associated random tree sequence.
Let $(T,r)$ be a numbered tree.
As $inj_{lab}((T,r),{\bm T}_n)$ is bounded, we have that $inj_{lab}((T,r),{\bm T}_n)$ is convergent for every $(T,r)$ if and only if we have that
\begin{equation}
\label{eqn.lab.conc}
\mathds{D}^2\left ({X_n^T\over n}\right )=\mathds{D}^2(inj_{lab}((T,r),{\bm T}_n))\rightarrow 0.
\end{equation}
We will use this formula to prove properties of the degree sequence.
We can expand the above formula using (\ref{eqn.sum.of.indicator})
\begin{multline}
\label{eqn.expand.deviation}
\mathds{D}^2\left({X_n^T\over n}\right)={1\over n^2}
\Bigl(\sum_{\psi,\phi: V(T)\mapsto [n]}\mathds{E}({\bm I}_n(T,\psi){\bm I}_n(T,\phi))-\\
-\sum_{\psi,\phi: V(T)\mapsto [n]}\mathds{E}({\bm I}_n(T,\psi))\mathds{E}({\bm I}_n(T,\phi))\Bigr)\rightarrow 0.
\end{multline}
The following two lemmas will establish a connection between the degree sequence and the probabilities
$\mathds{P}({\bm I}_n(T,\phi)=1)$.
Then we will use (\ref{eqn.expand.deviation}) to prove that the degree sequence satisfies the conditions in Lemma \ref{lem.exchange.degen}.
\begin{remark}
\label{rem.relab}
As the degree sequence is exchangeable we have that for any $\psi,\phi: V(T)\mapsto [n]$
$$\mathds{P}({\bm I}_n(T,\phi)=1)=\mathds{P}({\bm I}_n(T,\psi)=1).$$
\end{remark}
\noindent
Let $T$ be an arbitrary tree on $k$ nodes.
For a random degree sequence ${\bm D}_n$ and $\phi:V(T)\mapsto [n]$ let ${\bm D}_{\phi}=({\bm D}_n(\phi(i)))_{i=1}^k$.

\begin{lemma}
\label{lem.forest.prob}
Let ${\bm D}_n\in {\cal D}_n$ be a random degree sequence and ${\bm T}_n={\bm T}({\bm D}_n)$.
Let $F$ be an arbitrary forest on $m\ (m\leq n)$ nodes with remainder degrees $r=(r_1,\cdots , r_m)$.
Let $R=\sum_ir_i$ and denote by $C_1,C_2,\cdots ,C_c$ the connected components of $F$. The probability that an arbitrary $\phi:V(T)\mapsto [n]$ is an injective labeled homomorphism is
$$
\mathds{P}({\bm I}_n(F,\phi)=1)={(n-m+c-2)!\over (n-2)!}H(r,F)\mathds{P}({\bm D}_{\phi}=D_T),
$$
where $H(r,F)=\prod_{i=1}^c\left[\left(\sum_{j\in C_i}r_j\right)\prod_{j\in C_i}{(D_F(j)+r_j-1)!\over (r_j!)}\right]$ is a constant depending only on $F$ and the remainder degrees $r$.
\end{lemma}

\noindent
{\bf Proof:}
We may assume that $\phi(i)=i,\forall i\in V(F)$.
Let $R_i=\sum_{j\in C_i}r_j$.
Fix a degree sequence $D=(D(i))_{i=1}^n$.
It follows from the Pr\H{u}fer sequence that the number of trees realizing this degree sequence is ${n-2\choose D(1)-1,\cdots ,D(n)-1}$.
We need to count the trees with degree sequence $D$ which have $F$ spanned by the first $m$ nodes and the remainder degree condition holds.
Contract every connected component $C_i$ of $F$ to a single vertex $u_i$.
Also contract the images of these components in ${\bm T}({\bm D}_n)$.
We get a tree on $n-m+c$ nodes with degree sequence
$$D'=(R_1, R_2, \cdots , R_c,D(m+1),\cdots , D(n)).$$
There are
$$n-m+c \choose R_1-1, R_2-1, \cdots , R_c-1,D(m+1)-1,\cdots , D(n)-1$$
trees realizing the degree sequence $D'$.
For each connected component $C_i$ we can connect the $R_i$ edges to the vertices in
$$R_i!\over \prod_{j\in C_i}r_j!$$
ways.
It follows that the number of labeled trees realizing the degree sequence $D$ and having $F$ on the first $m$ vertices is
$$
{n-m+c-2\choose R_1-1,\cdots ,R_c-1,D(m+1)-1,\cdots ,D(n)-1}
\prod_{i=1}^c\left[ {R_i!\over \prod_{j\in C_i}r_j!}\right].
$$
From this it follows that
\begin{multline}
\label{eqn.forest.prob}
\mathds{P}({\bm I}_n(F,\phi)=1|{\bm D}_n=D)=\\
={\displaystyle {n-m+c-2\choose R_1-1,\cdots ,R_c-1,D(m+1)-1,\cdots ,D(n)-1}\over \displaystyle
{n-2\choose D(1)-1,\cdots ,D(n)-1}}\prod_{i=1}^c {R_i!\over \prod_{j\in C_i}r_j!}.
\end{multline}

Note that the degree sequence $D$ should be such that $D(i)=D_F(i)+r_i,\,i=1,\cdots,m$ holds for the first $m$ degrees.
We need to sum this probability for every possible degree sequence.
In our case we sum over degree sequences for which $D(i)=D_F(i)+r_i,\,i=1,\cdots,m$ holds.
As in equation (\ref{eqn.forest.prob}) the right hand side does not depend on $D(i),\,i>m$, we have
\begin{multline*}
\mathds{P}({\bm I}_n(F,\phi)=1)=\\
={(n-m+c-2)!\over (n-2)!}
\prod_{i=1}^c\left[R_i\prod_{j\in C_i}{(D_F(j)+r_j-1)!\over (r_j!)}\right]
\mathds{P}({\bm D}_{\phi}=D_T).
\end{multline*}
If we take $H(r,F)=\prod_{i=1}^c\left[R_i\prod_{j\in C_i}{(D_F(j)+r_j-1)!\over (r_j!)}\right]$, we get the desired equation.

\qed

Let $(F_1,r_1),(F_2,r_2)$ be two labeled graphs, $\phi:V(F_1)\mapsto [n]$ and $\psi:V(F_2)\mapsto [n]$.
We denote by $F_{1,2}$ the graph obtained by identifying nodes $i\in V(F_1),\ j\in V(F_2)$ if and only if $\phi(i)=\psi(j)$.
We can define remainder degrees $r_{1,2}$ on $F_{1,2}$ in a straightforward way if $\phi(i)=\psi(j)\Rightarrow r_1(i)=r_2(j)$.

\begin{lemma}
\label{lem.forest.cond.prob}
Let ${\bm D}_n\in {\cal D}_n$ be a random degree sequence and ${\bm T}_n={\bm T}({\bm D}_n)$.
Let $F_1,F_2$ be two forests on $m_1$ and $m_2$ nodes $(m_1,m_2\leq n)$ with remainder degrees $r_1,r_2$.
Let $\phi:V(F_1)\mapsto [n]$ and $\psi:V(F_2)\mapsto [n]$.
If $F_{1,2}$ is a forest and we can define $r_{1,2}$, then let
$m_{1,2}=|V(F_{1,2})|$, $c_{1,2}=\{$the number of components of $F_{1,2}\}$ and $R_{1,2}=\sum_{V(F_{1,2})}r_{1,2}(i)$. We have
\begin{multline}
\label{eqn.cond.forest.prob}
\mathds{P}({\bm I}_n(F_1,\phi)=1|{\bm I}_n(F_2,\psi)=1)=\\
{(n-m_{1,2}+c_{1,2}-2)! \over (n-m_2+c_2-2)!}
{H(r_{1,2},F_{1,2})\over H(r_2,F_2)}
\mathds{P}({\bm D}_{\phi}=D_{F_1}|{\bm D}_{\psi}=D_{F_2}).
\end{multline}

\end{lemma}
\noindent
{\bf Proof:}
The proof follows immediately from the definition of conditional probability.

\qed

Let ${\bm D}_n$ be a degree sequence and ${\bm T}_n={\bm T}({\bm D}_n)$ be the associated random tree.
Assume that $inj_{lab}((T,r),{\bm T}_n)$ is convergent.
Then by (\ref{eqn.lab.conc}) we have that 
$\mathds{D}^2({\bm X}_n^{T}\slash n)\to 0$.
For any tree $T$ on $k$ nodes we have
\begin{equation}
\label{eqn.2nd.moment}
\mathds{D}^2\left({{\bm X}_n^T\over n}\right)={1\over n^2}\sum_{\phi,\psi:V(T)\mapsto [n]}
\Big( \mathds{E}({\bm I}_n(T,\phi){\bm I}_n(T,\psi))-\mathds{E}({\bm I}_n(T,\phi))\mathds{E}({\bm I}_n(T,\psi))\Big)
\end{equation}
Now if we split the sum by the size of the intersection of $\phi(V(T))$ and $\psi(V(T))$ and use Remark \ref{rem.relab}, we have

\begin{multline}
\label{eqn.2nd.moment.b}
\mathds{D}^2\left({{\bm X}_n^T\over n}\right)=
{1\over n^2}\sum_{\substack{i=0\\|\phi(V(T))\cap \psi(V(T))|=i}}^k n(n-1)\cdot\ldots\cdot (n-2k+i+1)\cdot\\
\cdot\big(\mathds{E}({\bm I}_n(T,\phi){\bm I}_n(T,\psi))-\mathds{E}({\bm I}_n(T,\phi))\mathds{E}({\bm I}_n(T,\psi))\big)
\end{multline}
From Lemma \ref{lem.forest.prob} and \ref{lem.forest.cond.prob} we can easily derive that the order of the terms corresponding to $i\neq 0$ is ${\cal O}({1\over n})$.
It follows that the condition $\mathds{D}^2({\bm X}_n^{T})\to 0$ is equivalent to
$$
{(n-1)\cdot\ldots\cdot (n-2k+1)\over n}\left(\mathds{P}({\bm I}_n(T,\phi){\bm I}_n(T,\psi))-\mathds{P}({\bm I}_n(T,\phi))\mathds{P}({\bm I}_n(T,\psi))\right)\to 0.
$$
Using again Lemma \ref{lem.forest.prob} and \ref{lem.forest.cond.prob} we can easily derive the following:
\begin{multline}
\label{eqn.independency}
\forall T,\ \mathds{D}^2\left({{\bm X}_n^{T}\over n}\right)\to 0 \Leftrightarrow 
\forall \phi,\psi:V(T)\mapsto [n],\, \phi(V(T))\cap \psi(V(T))=\emptyset\\
\mathds{P}({\bm D}_{\phi}=D_{T}, {\bm D}_{\psi}=D_{T})\to \mathds{P}({\bm D}_{\phi}=D_{T})\mathds{P}({\bm D}_{\psi}=D_{T})
\end{multline}
\noindent
The following corollary is an easy application of Lemma \ref{lem.exchange.degen} and (\ref{eqn.independency}).

\begin{corollary}
\label{cor.conv.deFinetti}
The labeled subgraph densities of a random tree sequence converge in probability if and only if the corresponding degree sequence converges to an infinite IID sequence.
\end{corollary}
\begin{remark}
\label{rem.deg.conn}
The formula in Lemma \ref{lem.forest.prob} yields an easy result on the probability that two vertices $i,j$ with degrees $d_i,d_j$ are connected:
$$
\mathds{P}(ij\in E({\bm T}({\bm D}_n))\ |\ {\bm D}_n(i)=d_i,{\bm D}_n(j)=d_j)={d_i+d_j-2\over n-2}.
$$
Similarly for a given edge $ij\in E({\bm T}({\bm D}_n))$ the degree distribution of the vertices $i$ and $j$ can be expressed:
\begin{multline*}
\mathds{P}({\bm D}_n(i)=d_i,{\bm D}_n(j)=d_j\ |\ ij\in E({\bm T}({\bm D}_n)))=\\
{n\over n-2}{d_i+d_j-2\over 2}\mathds{P}({\bm D}_n(i)=d_i,{\bm D}_n(j)=d_j)
\end{multline*}
\end{remark}


\section{The limit of ${\bm T}({\bm D}_n)$}
\label{sec.limit}

In the last section we discussed tree sequences ${\bm T}_n$ for which $inj_{lab}((T,r),{\bm T}_n)$ was convergent for every $(T,r)$.
We now turn to neighborhood statistics.
First we want to relate them to labeled subgraph densities.
We will express the neighborhood statistics as functions of the labeled subgraph densities.

As before let $U^l$ denote the set of all finite $l$-deep rooted tree.
Consider an $l$-deep rooted tree with root $x$: $T^l_x\in U^l$, with $|T^l_x|=k$.
Let us denote the nodes at distance $i$ from the root by $T_i$, and $|T_i|=t_i$ ($t_0$ is just $1$, $t_1$ is the degree of the root).
$B_G(v,l)$ is the rooted $l$-ball around $v$ in $G$ and ${\bm T}_n={\bm T}({\bm D}_n)$ is a random labeled tree with degree distribution ${\bm D}_n$.

Let $\sigma, \rho\in Aut(T_x^l)$ be two rooted automorphisms of the rooted tree $T_x^l$.
We say that $\sigma \sim \rho$ if and only if there exists $\tau \in Aut(T_x^l)$, such that $\tau$ fixes every vertex not in $T_{l}$ and $\sigma \circ \tau = \rho$.
$\sim$ is an equivalence relation.
The equivalence classes have $\prod_{i\in T_{l-1}}(D(i)-1)!$ elements, hence it follows
\begin{equation}
\label{eqn.aut}
|Aut(T_x^l)|=|Aut(T_x^l)\slash\sim|\prod_{i\in T_{l}}(D(i)-1)!.
\end{equation}
It is easy to see that
\begin{equation}
\label{eqn.neigh.lab.equiv}
p(l,T_x^l,{\bm T}_n)={1\over n}{\displaystyle {\bm X}_n^{(T',(r_i')_{i=1}^{|T'|})}\over |Aut(T_x^l)\slash \sim|}, \textrm{ where}
\end{equation}
\begin{equation}
\label{eqn.neig.lab.eqv}
\begin{array}{l}
T'=T_x^l\setminus T_l \\
r_i'=\left \{
\begin{array}{ll}
0 &  i\notin T_l\cup T_{l-1}\\
D_{T_x^l}(i)-1 & i\in T_{l-1}.
\end{array}
\right.
\end{array}
\end{equation}

If ${\bm T}_n$ is a convergent random tree sequence then from (\ref{eqn.rnd.conv}) and (\ref{eqn.neigh.lab.equiv}) we have that for any $T'$ defined above
\begin{equation}
\mathds{D}^2\left({X_n^{T'} \over n}\right) \to 0.
\end{equation}
%
%
For bounded degree graphs the convergence of the neighborhood densities implies the convergence of the graph sequence in the sense of Benjamini and Schramm.
We saw earlier in Example \ref{example.star} that for bounded average degree graphs this is not the case.
The convergence of the neighborhood densities alone is not enough.
We need also (\ref{eqn.neigh.cons}) to hold.

The reason is that for fixed $k$ the $k$-neighborhood of the large degree nodes is large (${\cal O}(n)$).
In Example \ref{example.star} even if $k=1$, every node "sees" the center node (eg. every neighborhood with radius $1$ contains the center node) and so every $2$ radius neighborhood contains ${\cal O}(n)$ vertices, which is unbounded.

Assign remainder degrees $r$ ($r_i=0,\ \forall i\notin T_l$) to the rooted tree $T_x^l$ and forget the root, then using Lemma \ref{lem.forest.prob}
\begin{multline}
\label{eqn.exp.subtree}
\mathds{E}({\bm X}_n^T)=\mathds{E}\left(\sum_{\phi:V(T)\mapsto [n]}{\bm I}_n(T,\phi)\right)={n!\over (n-k)!}\mathds{P}({\bm I}_n(T,\phi)=1)=\\
n{n-1\over n-k}\mathds{P}({\bm D}_n(\{1,2,\cdots ,k\})=D_T)H(r,T).
\end{multline}
From (\ref{eqn.neig.lab.eqv}) we have that
$$
p(l,T_x^l,{\bm T}_n)={1\over n}{{\bm X}_n^{T'}\over |Aut(T_x^l)\slash\sim|}.
$$
We want to define an infinite random rooted tree which is the limit of ${\bm T}_n$.
Let
$$
\mu_n(T_x^l)={1\over n}{\mathds{E}({\bm X}_n^{T'})\over |Aut(T_x^l)\slash\sim|}.
$$

Assume we have a convergent sequence of random trees ${\bm T}_n$ with degree sequence ${\bm D}_n$.
Further assume that ${\bm D}_n\to {\bm D}=({\bm D}_0,{\bm D}_0,\cdots)$ and let $\gamma=\mathds{E}(D_0)-1$.
Define

\begin{multline}
\label{eqn.prob.dist}
p(T_x^l)=\lim_{n\to\infty}\mu_n(T_x^l)=\\
\lim_{n\to\infty}{1\over |Aut(T_x^l)\slash\sim|}
{n-1\over n-k}\mathds{P}({\bm D}_n(\{1,2,\cdots ,k\})=D_{T'})H(r,T')=\\
{\prod_{i\notin T_{l}}\mathds{P}(D_0=d_i)(d_i-1)!\over |Aut(T_x^{l})|}t_l
\end{multline}
We can expand the formula
\begin{multline*}
H(r,T')=\sum_{i\in V(T')}r_i\prod_{i\in V(T')}{(d_i+r_i-1)!\over r_i!}=\\
\sum_{i\in T_{l-1}}(d_i-1)\prod_{i\notin T_{l-1}\cup T_l}(d_i-1)!=t_l\prod_{i\notin T_{l-1}\cup T_l}(d_i-1)!.
\end{multline*}
Then the last equation in (\ref{eqn.prob.dist}) follows using equation (\ref{eqn.aut}) and the expansion of $H(r,T')$.

Define $\mu({\cal T}(F))=p(F)$.
As the sets ${\cal T}(F)$ generate the $\sigma$-algebra, we can extend $\mu$ to $\cal T$ if $\mu$ satisfies (\ref{eqn.neigh.cons}).
If this is the case then $\mu$ is a random infinite rooted tree.
\begin{lemma}
\label{lem.mu.prob.measure}
Let ${\bm D}_n$ be an exchangeable random degree sequence and assume that ${\bm D}_n\to {\bm D}$, where ${\bm D}$ is an infinite IID random sequence of the variable ${\bm D}_0$.
Let $\mu$ be the associated measure defined above.
$\mu$ extends to a probability measure on $\cal G$ if and only if
$\mathds{E}({\bm D}_0)=2$ (or equivalently $\gamma=1$).
\end{lemma}

\noindent
{\bf Proof:}
We only need to show that $\mu$ satisfies (\ref{eqn.neigh.cons})
\begin{equation}
\label{eqn.prob.measure}
p(T_x^{l-1})=\sum_{\displaystyle T_x^{l}: B_{T_x^{l}}(x,l-1)\cong T_x^{l-1}}p(T_x^{l})\Leftrightarrow \gamma=1
\end{equation}
We have
$$
p(T_x^l)={\prod_{i\notin T_l}\mathds{P}({\bm D}_0=d_i)(d_i-1)!\over |Aut(T^l_x)|}t_l.
$$
Now rearranging the sum by the degrees of the leafs of $T_x^{l-1}$ in $T_x^l$ we have

\begin{multline*}
\sum_{\displaystyle T_x^{l}: B_{T_x^{l}}(x,l-1)\cong T_x^{l-1}}p(T_x^{l})=
\sum_{D_{T_x^l}(i)=1,\,i\in T_{l-1}}^\infty p(T_x^{l-1}\cup (d_i)_{i\in T_{l-1}})=\\
\prod_{j\notin T_{l-1}\cup T_l}\mathds{P}({\bm D}_0(j)=d_j)(d_j-1)!
\sum_{D_{T_x^l}(i)=1,\,i\in T_{l-1}}^\infty {\prod_{i\in T_{l-1}}\mathds{P}({\bm D}_0(j)=d_j)(d_i-1)!\over |Aut(T_x^l)|}t_l=\\
{\prod_{j\notin T_l\cup T_{l-1}}\mathds{P}({\bm D}_0(j)=d_j)(d_j-1)!\over
|Aut(T_x^l\setminus T_l)|}t_{l-1}\gamma,
\end{multline*}
where the last equation follows from the fact that for fixed $d_i,\ i\in T_{l-1}$ every $\sigma\in Aut(T_x^{l-1})$ has only one extension in $Aut(T_{l})\slash\sim$.
Now (\ref{eqn.prob.measure}) will hold only if $\gamma=1$.
It follows that (\ref{eqn.prob.measure}) holds if and only if $\mathds{E}(D_0)=2$ $(\gamma=1)$.

\qed

\noindent
{\bf Proof of Theorem \ref{thm.main}:}\\
\noindent
Let ${\bm D}_n$ be a degree sequence and ${\bm T}({\bm D}_n)={\bm T}_n$ be the associated random tree sequence.
First assume, that the degree sequence converges to the distribution ${\bm D}=({\bm D}_0,{\bm D}_0,\cdots)$ and $\mathds{E}({\bm D}_0)=2$.
From equation (\ref{eqn.independency}) we get that for an arbitrary tree $T$, $\mathds{D}^2\left({{\bm X}_n^T\over n}\right)\rightarrow 0$.
Then by equation \eqref{eqn.neigh.lab.equiv} we have that for every $T_x^l$ $l$-deep rooted tree, the neighborhood statistics converge in probability to a limiting distribution $p(T_x^{l})$.
As the assumptions of Lemma \ref{lem.mu.prob.measure} hold we have that $p(T_x^l)$ defines a measure $\mu$ on infinite rooted trees and so ${\bm T}_n\rightarrow \mu$.

On the other hand assume that ${\bm T}_n$ converges to a random infinite rooted tree $\mu$.
Then by equation (\ref{eqn.neigh.lab.equiv}) we get that the number of degree $d$ vertices is concentrated.
Using that our degree distribution is exchangeable we get that ${\bm D}_n\rightarrow {\bm D}=({\bm D}_0,{\bm D}_0,\cdots )$ and $\mathds{E}({\bm D}_0)=2$.
This completes the proof of Theorem \ref{thm.main}.

\end{document}